\documentclass[11pt,reqno]{amsart}
\usepackage[matrix,arrow,tips,curve]{xy}
\usepackage[english]{babel}
\usepackage{amssymb,amsthm}
\usepackage{enumerate}
\usepackage{amsmath} 

\newtheorem{thm}{Theorem}
\newtheorem{lemma}{Lemma}
\newtheorem{prop}{Proposition}
\newtheorem{cor}{Corollary}

\theoremstyle{definition}
\newtheorem{defi}{Definition}

\theoremstyle{remark}
\newtheorem{remark}{Remark}
\newtheorem{example}{Example}

\newcommand{\ov}{\overline}

\newcommand{\ha}{\widehat}
\newcommand{\ma}{\mathcal}

\newcommand{\p}{\mathbb{P}}
\newcommand{\ol}{\mathcal{O}}

\newcommand{\Z}{\mathbb{Z}}
\newcommand{\Q}{\mathbb{Q}}

\newcommand{\R}{\mathbb{R}}

\newcommand{\Supp}{\operatorname{Supp}}
\newcommand{\Sing}{\operatorname{Sing}}
\newcommand{\Pic}{\operatorname{Pic}}

\newcommand{\im}{\operatorname{Im}}
\newcommand{\NE}{\operatorname{NE}}

\newcommand{\Chow}{\operatorname{Chow}}
\newcommand{\Lo}{\operatorname{Locus}}

\renewcommand{\thesubsection}{\arabic{subsection}}

\begin{document}
\title{On covering and quasi-unsplit
families of rational curves}
\author[L.\ Bonavero, C.\ Casagrande, S.\ Druel]{Laurent BONAVERO,
  Cinzia CASAGRANDE, St\'ephane DRUEL} 
\date{1st April 2005}
\maketitle
\noindent
\def\restriction{\string |}

\newcommand\mero{{\dashrightarrow}}
\def\thesubsub{\thesubsection .\arabic{subsub}}
\def\subsub#1{\addtocounter{subsub}{1}\par\vspace{3mm}
\noindent{\bf \thesubsub ~ #1 }\par\vspace{2mm}}
\def\coker{\mathop{\rm coker}\nolimits}
\def\pr{\mathop{\rm pr}\nolimits}
\def\im{\mathop{\rm Im}\nolimits}
\def\hfl#1#2{\smash{\mathop{\hbox to 12mm{\rightarrowfill}}
\limits^{\scriptstyle#1}_{\scriptstyle#2}}}
\def\vfl#1#2{\llap{$\scriptstyle #1$}\big\downarrow
\big\uparrow
\rlap{$\scriptstyle #2$}}
\def\diagram#1{\def\normalbaselines{\baselineskip=0pt
\lineskip=10pt\lineskiplimit=1pt}   \matrix{#1}}
\def\limind{\mathop{\oalign{lim\cr
\hidewidth$\longrightarrow$\hidewidth\cr}}}

\long\def\InsertFig#1 #2 #3 #4\EndFig{
\hbox{\hskip #1 mm$\vbox to #2 mm{\vfil\includegraphics{#3}}#4$}}
\long\def\LabelTeX#1 #2 #3\ELTX{\rlap{\kern#1mm\raise#2mm\hbox{#3}}}

{\let\thefootnote\relax
\footnote{%\hskip3em 
\indent \textbf{A.M.S.~classification :} 14E30, 14J99, 14M99. 
}}

\vspace{-1cm}

\begin{center}
\begin{minipage}{130mm}
\scriptsize

{\bf Abstract.} 
We study extremality properties of covering 
families of rational curves on projective varieties. Among others,
we show that on a normal and $\Q$-factorial projective variety
$X$ with $\dim(X) \leq 4$, every covering and quasi-unsplit
family of rational curves generates a geometric 
extremal ray of the Mori cone 
$\overline{\rm NE}(X)$ of classes of effective 1-cycles. 

\end{minipage}
\end{center}

\section{Introduction}

\noindent Let $X$ be a normal and uniruled
complex projective
variety. 
Consider an irreducible and closed subset
$V$ of $\Chow(X)$ such that:
\begin{enumerate}[$\circ$]
\item any element of $V$ is a cycle whose irreducible components are
rational curves; 
\item $V$ is covering (which means that for any point
$x \in X$, there exists an element of $V$ passing through $x$).
\end{enumerate}
We call such a $V$ a \emph{covering family of rational 1-cycles} on $X$.
If moreover,
all irreducible
components of the cycles parametrized by $V$ 
are numerically proportional, 
 we call $V$ a 
\emph{covering and quasi-unsplit 
family of rational 1-cycles} on $X$ 
(see \cite[Definition 2.13]{choc}).

For any covering family $V$ 
of rational 1-cycles on $X$, we will denote by
$[V]$ the numerical class in ${\rm NE}(X)$
of the general cycle
of the family $V$ and by $\R_{\geq 0} [V]$ the half-line
generated by $[V]$.

A \emph{geometric extremal ray} of the Mori cone
$\overline{\NE}(X)$ is
a half-line $R\subseteq \overline{\NE}(X)$ such that if
$\gamma_1+\gamma_2\in R$ for some
$\gamma_1,\gamma_2\in \overline{\NE}(X)$, then $\gamma_1,\gamma_2\in
R$
(see Section 2 for precise definitions and notation).

\medskip

\noindent{\bf Question.}
Let $V$ be a covering and quasi-unsplit family of rational 1-cycles on
$X$. Is $\R_{\geq 0} [V]$ a geometric extremal ray
of $\overline{\rm NE}(X)$?

\medskip

Note that this question is natural, since any
family of rational 1-cycles such that the general member 
generates a geometric extremal ray of $\overline{\NE}(X)$
is quasi-unsplit. The converse is not true if the family
is not covering (just think of a smooth blow-down of
a smooth projective variety to a non projective one).

\medskip  

Let $V$ be any covering family of rational 1-cycles on $X$. 
Then $V$ defines set-theoretically an equivalence
relation on $X$: two points $x,x'$ are 
\emph{$V$-equivalent} if there exist $v_1,\dotsc,v_m\in V$ such that
some connected component of $C_{v_1}\cup\cdots\cup C_{v_m}$
contains $x$ and $x'$, where $C_v\subset X$ is the curve
 corresponding to $v\in V$.

In this situation, after Campana's results (see Section 2),
there exists an almost holomorphic map $q\colon X\dasharrow Y$,
to a projective algebraic variety, 
whose general fibers are $V$-equivalence classes.

We first prove the following result,
which involves the dimension
of the general fiber of $q$.
\begin{thm}
\label{result}
Let $X$ be a normal and $\Q$-factorial complex projective
variety of dimension $n$.
Let $V$ be  a
covering and quasi-unsplit family of rational 1-cycles on $X$, and 
let $f_V$ be the dimension of a general $V$-equivalence class.

If $f_V \geq n-3$, then $\R_{\geq 0} [V]$ is a geometric extremal ray
of the Mori cone $\overline{\rm NE}(X)$.
\end{thm}
We then immediately get the following.
\begin{cor}
Let $X$ be a normal and $\Q$-factorial complex projective
variety of dimension $n\leq 4$.
Let $V$ be a
covering and quasi-unsplit family of rational 1-cycles on $X$.
Then $\R_{\geq 0} [V]$ is a geometric extremal ray
of the Mori cone $\overline{\rm NE}(X)$.
\end{cor}

\smallskip

As previously recalled, one can associate a rational map
$q\colon X\dasharrow Y$
to any covering family of rational 1-cycles on $X$. 
We call a \emph{geometric quotient} for $V$ a morphism
$q'\colon X\to Y'$, onto a normal projective variety $Y'$, 
such that \emph{every} fiber of $q'$ is a $V$-equivalence class.
If such a quotient exists, then it is clearly unique up to isomorphism. On
the other hand, even if $X$ is smooth,
a geometric quotient for $V$ does not necessarily exist 
(see example~\ref{Z2}).

\medskip

The study of the extremal contraction given by the previous result
leads to the following.
\begin{thm}
\label{toto}
Let $X$ be a normal and $\Q$-factorial complex projective
variety, having canonical singularities, of dimension $n$.
Let $V$ be  a
covering and quasi-unsplit family of rational 1-cycles on $X$, and 
let $f_V$ be the dimension of a general $V$-equivalence class.

If $f_V \geq n-3$, 
then the Mori contraction  of $\R_{\geq 0} [V]$,
${\rm cont}_{[V]}\colon X\to Y'$,  
is the geometric quotient for $V$. 
If moreover $f_V\geq n-2$, then ${\rm cont}_{[V]}$ is equidimensional.
\end{thm}
We finally consider the toric case, where we can prove both
extremality and existence
of the geometric quotient for a quasi-unsplit family in any dimension.
\begin{thm}
\label{toric}
Let  $X$ be a toric and $\Q$-factorial complex projective variety, and 
let $V$ be a
quasi-unsplit covering family of rational 1-cycles in $X$. 
Then $\R_{\geq 0} [V]$ is a geometric extremal ray
of the Mori cone $\overline{\rm NE}(X)$,
and the Mori contraction of $\R_{\geq 0} [V]$, 
${\rm cont}_{[V]}\colon X\to Y'$, 
is the geometric quotient for $V$.
\end{thm}
The following is an immediate application of Theorems \ref{result} and
\ref{toric}. 
\begin{cor}
Let $X\subset\p^N$ be a normal and $\Q$-factorial variety, 
covered by
lines. Assume either that $X$ is toric, or that 
$X$ has canonical singularities
and $\dim X\leq 4$. Let $V$ be an irreducible family of lines
covering $X$. 

Then there exists a morphism $q'\colon X\to Y'$, onto a normal,
$\Q$-factorial, 
projective variety $Y'$ with $\rho_{Y'}=\rho_X-1$,\footnote
{We denote by $\rho_Z$ 
the Picard number of an algebraic variety $Z$.}
such that all lines of $V$ are contracted by~$q'$.
\end{cor}
%%%%%%%%%%%%%%%%%%%%%%%%%%%%%%%%%%%%%%%%%%%%%%%%%%%%%%%%%%%%%%%%%%%%%%%%%%%%%%%
\section{Set-up on families of rational 1-cycles}
\noindent Let $X$ be a normal, irreducible, $n$-dimensional complex projective
variety. We denote by $\mathcal{N}_1(X)_{\R}$ (respectively,
$\mathcal{N}_1(X)_{\Q}$) the vector space of 1-cycles in $X$ with real
(respectively, rational) coefficients, modulo numerical equivalence.
In $\mathcal{N}_1(X)_{\R}$, let $\overline{\NE}(X)$ be the closure of
the cone generated by classes of effective 1-cycles in $X$. 

Recall that the existence of a covering 
family $V$ of rational 1-cycles on $X$ is equivalent to $X$ being
uniruled \cite[Proposition IV.1.3]{kollar}.

For such family $V$, we
have a diagram given by the incidence variety
$\ma{C}$ associated to~$V$:
\begin{equation}\label{D1}
\xymatrix{ {\ma{C}}\ar[d]^{\pi}\ar[r]^F & X \\ V & } 
\end{equation}
where $\pi$ and $F$ are proper and surjective. 
We set $C_v:=F(\pi^{-1}(v))$ for any $v\in V$.
 
The relation of $V$-equivalence on $X$ induced by such a family was
introduced and studied in \cite{campana81}; we refer the reader to
\cite{campa},
\cite[\S 5.4]{debarreUT} or \cite[\S IV.4]{kollar} for more details.
In particular, there exists 
a rational map $q\colon X\dasharrow Y$ associated
to $V$, whose main properties are recalled now. By 
\cite[Theorem 5.9]{debarreUT},  
there exists
a closed and irreducible subset of $\Chow (X)$
whose normalization $Y$ 
satisfies the following properties:
\begin{enumerate}[(1)]
\item[(a)] let ${Z}\subset Y\times X$ be
the restriction of the universal family, 
\begin{equation}\label{D2}
 \xymatrix{ {{Z}}\ar[d]^{p}\ar[r]^e & X \ar@{-->}[ld]^q
\\ Y & } \end{equation}
then $e$ is birational and $q = p \circ e^{-1}$ is almost holomorphic
(which means that the indeterminacy locus of $q$ 
does not dominate $Y$);
\item[(b)] a general fiber of $q$ is a $V$-equivalence class, 
\item[(c)] a general fiber of $q$, hence of $p$, is irreducible.
\end{enumerate}

As a consequence of the existence of this map $q$, 
a general $V$-equivalence class is a closed
subset of $X$. We denote by $f_V$ its dimension, so that $\dim
Y=n-f_V$. 
Moreover, it is well known that any $V$-equivalence
class is a countable union of
closed subsets of $X$.  
\begin{defi}
We say that a subset $Z$ of $X$ is \emph{$V$-rationally connected} if
every connected component of $Z$ is contained in some
$V$-equivalence class. 
\end{defi}
\begin{lemma} 
\label{ultimo}
Let $X$ be a normal projective
variety and 
$V$ be a covering 
family of rational 1-cycles on $X$.
Consider the diagram \eqref{D2} above.
Then $e(p^{-1}(y))$ is $V$-rationally connected for any $y\in Y$.
\end{lemma}
\begin{proof} Let $\mathcal R \subset X\times X$ the graph
of the equivalence relation defined by $V$: it is a countable
union of closed subvarieties since $V$ is proper.  
The fiber product $Z\times_{Y} Z$ is irreducible and
thus $(e\times e) (Z\times_{Y} Z) \subset \mathcal R$
thanks to properties (a) and (b) above. 
Therefore, for any $x \in e(p^{-1}(y))$, the cycle 
$e(p^{-1}(y))$ is contained in the $V$-equivalence class
of $x$.
\end{proof}
The following well known remark will be of
constant use (see \cite[Proposition IV.3.13.3]{kollar}, 
or \cite[Corollary 4.2]{occhettaGM}).
\begin{remark}
\label{rc}
If $Z \subset X$ is $V$-rationally connected, every curve 
contained in $Z$ is
numerically equivalent in $X$ to a linear combination with rational
coefficients of irreducible components of cycles in $V$.
In particular, if $V$ is quasi-unsplit, 
the numerical class of every curve contained in
a $V$-rationally connected subset $Z$ of $X$ belongs to 
$\R_{\geq 0} [V]$.
\end{remark}
Finally, we will need the following.
\begin{lemma} 
\label{families}
Let $X$ be a normal projective
variety and 
$V$ be a covering and quasi-unsplit
family of rational 1-cycles on $X$. 
Then there exists a covering and quasi-unsplit
family ${V}'$ of rational 1-cycles on $X$ such that:
\begin{enumerate}[$\circ$]
\item the general cycle of ${V}'$ is reduced and irreducible;
\item for any
${v}'\in{V}'$ there exists $v\in V$ such that
$C_{{v'}}\subseteq 
C_v$; in particular $\R_{\geq 0} [V]= \R_{\geq 0} [{V'}]$.
\end{enumerate}
\end{lemma}
\begin{proof}
Let $\ma{C}$ be the incidence variety
associated to~$V$ as in \eqref{D1}.

It is well-known that every irreducible component
of ${\ma{C}}$ dominates $V$, 
let $\ma{C}''$ be an irreducible component
of ${\ma{C}}$ which dominates $X$ too.
Let ${\ma{C}'}$ be the normalization of 
$\ma{C}''$ and $\ma{C}' \to {V}'$ be the Stein
factorization 
of the composite map $\ma{C}'\to \ma{C}'' \to V$.
Since $\ma{C}' \to {V}'$ has connected fibers and
${\ma{C}'}$ is normal, the general fiber
of ${\ma{C}'} \to {V}'$ is irreducible.
Moreover, the image in $X$ of every fiber of   $\ma{C}' \to {V}'$  is
contained in a cycle of $V$.

Since ${V}'$ is normal, there is a holomorphic map
${V}' \to {\rm Chow}(X)$. Then after replacing ${V}'$ by its
image in ${\rm Chow}(X)$ and 
${\ma{C}'}$ by its image in $\Chow(X)\times X$, 
we get the desired family.  
\end{proof}
%%%%%%%%%%%%%%%%%%%%%%%%%%%%%%%%%%%%%%%%%%%%%%%%%%%%%%%%%%%%%%%%%%%%%%%
\section{Properties of the base locus and extremality}
\noindent Let $V$ be a covering family of rational 1-cycles on $X$, and
recall the diagram \eqref{D2} associated to $V$.

Let $E\subset{Z}$ be the exceptional locus of $e$, and $B:=e(E)\subset 
X$. Observe that since $X$ is normal, $\dim B\leq n-2$. 
\begin{prop} 
\label{baselocus}
Let $X$ be a normal and $\Q$-factorial projective variety, and $V$ be
a covering and quasi-unsplit family of rational 1-cycles on $X$.
Consider the associated diagram as in \eqref{D2}.
Then:
\begin{enumerate}[$(i)$]
\item $e(p^{-1}(y))$ is a
 $V$-equivalence class of dimension $f_V$  for every $y\in Y\setminus p(E)$; 
\item $B$ is the union of all
$V$-equivalence classes of dimension bigger than $f_V$.
\end{enumerate}
\end{prop}
\begin{proof}
Set $X^0 : = X \setminus B$ and $Y^0: = Y \setminus p(E) = q(X^0)$.
Choose a very ample line bundle $L$ on $Y$, and 
let $U\subset |L|$ be the open subset of divisors $H$ that are
irreducible and such that $H \cap Y^0\neq \emptyset$.
For any $H$ in $U$,
we define $\ha{H}:=\ov{q^{-1}(H\cap Y^0)}$, which is a Weil
divisor in $X$. Since $X$ is $\Q$-factorial, 
some multiple of $\ha{H}$ defines a line bundle
$\ha{L}$ on $X$.

Let now $N:=h^0(L)$, 
and let $s_1,\dotsc,s_N$
be general global sections generating $L$. For each $i=1,\dotsc,N$, let
$H_i\in |L|$ be the divisor of zeros of $s_i$
and $\widehat{H}_i$ in $X$ as defined above.

Let's show that
$\widehat{H}_1\cap\cdots\cap\widehat{H}_N=B$. 
If $x\not\in B$, then ${q}$ is defined in $x$ and 
there is some $i_0\in\{1,\dotsc,N\}$ such that  ${q}(x)\not\in
H_{i_0}$, so $x\not\in\widehat{H}_{i_0}$. 
Conversely, let $x\in B$ and fix $i\in \{1,\dotsc,N\}$. 
Then $e^{-1}(x)$ has positive dimension; let $C\subset Z$
be an irreducible curve such that $e(C)=x$. Then $p(C)$ is a curve in
$Y$, hence 
$H_i\cap p(C)\neq \emptyset$ and $p^{-1}(H_i)\cap
C\neq\emptyset$. Now observe that $p^{-1}(H_i)$ does not contain any
component of $E$, hence $e(p^{-1}(H_i))$ is a divisor in $X$ which
coincides with  $\widehat{H}_i$ over $X\setminus B$. Then
$\widehat{H}_i=e(p^{-1}(H_i))$ and
$x\in \widehat{H}_i$.

Let $i\in\{1,\dotsc,N\}$.
Observe that $\widehat{H}_i\cdot [V]=0$, because $[V]$ is
quasi-unsplit and any irreducible
component of  
general cycle of the family is contained in a fiber of $q$
disjoint from  $\widehat{H}_i$.
This implies that  $\widehat{H}_i$  is closed with respect to
$V$-equivalence. In fact, let $C$ be an irreducible component of a
cycle
of $V$ such that
$C\cap \widehat{H}_i\neq\emptyset$. Since $V$ is quasi-unsplit,
we have
$\widehat{H}_i\cdot C=0$, which implies
$C\subseteq  \widehat{H}_i$.

Now since ${B}=\widehat{H}_1\cap\cdots\cap\widehat{H}_N$ and all 
$\widehat{H}_i$'s are closed with respect to $V$-equivalence, we see
that $B$ is a union of $V$-equivalence classes.

\medskip

Observe that if $C\subset X\setminus B$ is an irreducible curve such
that $\widehat{H}\cdot C=0$ for some $H\in U$, then $q(C)$ is a point. 
In fact,  if $q(C)$ is a curve,
there exists $H_0\in U$ such that $H_0$ intersects $q(C^0)$ in a
finite number of points.
Then $\ha{H}_0$ intersects $C$ without
containing it, a contradiction, because $\ha{H}$ and $\ha{H}_0$ are 
numerically equivalent, so $C\cdot \ha{H} > 0$.

Now fix $y_0\in Y^0$. We know by Lemma
\ref{ultimo} that $e(p^{-1}(y_0))$ is contained in a $V$-equivalence
class $F$. Since $B$ is closed with respect to $V$-equivalence, we
have $F\subset X^0$.
Consider an irreducible component $C$ of a cycle of $V$ 
such that $C\subseteq F$.
Since $V$ is quasi-unsplit, we have
$\widehat{H}\cdot C=0$, 
hence $q(C)$ is a point by what we proved above. 
Therefore
$q(F)=y_0$ and $F=e(p^{-1}(y_0))$, so we have $(i)$.

\medskip

For any $x\in X$, let $Y_x:=p(e^{-1}(x))$ be the family of 
cycles parametrized by $Y$ and
passing through $x$, and $\Lo(Y_x):=e(p^{-1}(Y_x))$. 
Observe that  for any $y\in Y_x$, the subset 
$e(p^{-1}(y))$ contains $x$ and is $V$-rationally connected by Lemma
\ref{ultimo}.
Hence
$\Lo(Y_x)$ is $V$-rationally connected for any $x\in X$. 

Since $Z\subset X\times Y$, we have $\dim Y_x=\dim e^{-1}(x)$. Thus
$\dim Y_x>0$ if and only if $x\in B$, 
by Zariski's main Theorem. If so, $\Lo(Y_x)$
has dimension at least $f_V+1$.

Now let $F$ be a $V$-equivalence class contained in $B$, and
$x\in F$. 
Then $\Lo(Y_x)$ has dimension at least $f_V+1$ and is contained in
$F$, hence $\dim F\geq f_V+1$.
\end{proof}
Let us  remark that in general, if
$V$ is not quasi-unsplit, $B$ is not closed with respect to $V$-equivalence.
\begin{example}
\label{Z2}
In $\mathbb P ^{2}$ fix two points $x$, $y$ and the line
$L=\ov{xy}$. 
Consider
$\mathbb P^{2}\times\mathbb P^{2}$ 
with the projections $\pi_1$, $\pi_2$ on the two
factors,  and fix
three curves $R_x$, $R_y$, $L'$ such that:
\begin{enumerate}[$\circ$]
\item $R_x$ is a line in $\mathbb P^{2}\times x$ and
$R_y$ is a line in $\mathbb P^{2}\times y$;
\item
$\pi_1(R_x)\cap\pi_1(R_y)$ is a point $z\in\mathbb P^{2}$;
\item
$L':= z\times L$ is the unique line dominating $L$ via $\pi_2$
and intersecting
both $R_x$ and $R_y$.
\end{enumerate}
Let $\sigma\colon W\to \mathbb P^{2}\times\mathbb P^{2}$ 
be the blow-up of $R_x$ and
 $R_y$. In 
$W$, the strict transform of $L'$ is a smooth rational curve
with normal bundle
$\ol_{\mathbb P^1}(-1)^{\oplus 3}$. Let $X$ be the variety obtained by
``flipping'' this curve.
Then $X$ is a smooth toric Fano $4$-fold with $\rho_X=4$ (this is $Z_2$ in
Batyrev's list, see \cite[Proposition 3.3.5]{bat2}).
$$ \xymatrix{
X \ar@{-->}[r]\ar@{-->}[rd]_q & W\ar[d]^{\pi_2\circ\sigma} \\
& {\mathbb P^{2}}}$$
The strict transform of a general line in a fiber of $\pi_2$ gives a
covering family $V$ of rational curves on $X$. The birational map
$X\dasharrow\p^2\times\p^2$ is an isomorphism over
$\p^2\times(\p^2\setminus L)$; if $U\subset X$ is the corresponding
open subset, then $U$ is closed with respect to $V$-equivalence and
every fiber of $q\colon U\to\p^2\setminus L$ is a
$V$-equivalence class isomorphic to $\p^2$. 
Thus $f_V=2$.

Let $T_x$ and $T_y$
be the images in $X$ of the exceptional divisors of
$\sigma$ in $W$. These two divisors are $V$-rationally connected, and
they can not be contained in $B$ because $\dim B\leq 2$. Moreover, 
$P:=T_x\cap T_y$ is the $\mathbb P^2$ with normal bundle $\ol_{\mathbb
  P^2}(-1)^{\oplus 2}$ obtained under the flip. 
The map $q\colon X\dasharrow\mathbb P^2$ can not be defined over $P$, so
$P\cap B\neq\emptyset$. Therefore
$B$ can not be closed with respect to $V$-equivalence.

Observe that the numerical class of $V$ lies in the interior of
$\overline{\NE}(X)$, hence the unique morphism, onto a
projective variety,  which contracts curves in $V$, is $X\to\{pt\}$.
\end{example}
A key observation is the following.
\begin{prop}
\label{extremal}
Let $X$ be a normal and $\Q$-factorial projective variety, and $V$ a
covering and quasi-unsplit family of rational 1-cycles on $X$.

If  $B$ is $V$-rationally connected, 
then $\R_{\geq 0} [V]$ is a geometric extremal ray of $\overline{\NE}(X)$.
\end{prop}
\begin{proof}
Let $X^0 : = X \setminus B$ and $Y^0 := Y \setminus p(E) = q(X^0)$.
Let $L$ be a very ample line bundle on $Y$. 
Let $U\subset |L|$ be the open subset of divisors $H$ that are
irreducible and such that $H \cap Y^0\neq \emptyset$.
For any $H$ in $U$,
we define $\ha{H}:=\ov{q^{-1}(H\cap Y^0)}$
as in the proof of Proposition \ref{baselocus}. 
Recall that $\ha{H}\cdot [V]=0$.

Let's show that $\ha{H}$ is nef. Assume by contradiction that there exists
an irreducible curve $C$ with $C\cdot \ha{H}<0$. 

\medskip 

\noindent {\em Claim.} $C\subseteq B$.

Actually, either $C$ is contained in a
fiber of $q$, hence it is numerically proportional to $[V]$
which contradicts $C\cdot \ha{H}<0$. Or 
$C\cap X^0 =: C^0$ is an open subset of $C$, 
$\dim q(C^0)=1$, hence
there exists $H_0\in U$ such that $H_0$ intersects $q(C^0)$ in a
finite number of points.
Then $\ha{H}_0$ intersects $C$ without
containing it, a contradiction, because $\ha{H}$ and $\ha{H}_0$ are 
numerically equivalent, so $C\cdot \ha{H} > 0$.

\medskip

Since  $B$ is $V$-rationally connected, 
$C$ must be  numerically proportional to $V$, impossible.

Let's finally show that $C\cdot\ha{H}=0$ if and only if $C$ is numerically 
proportional to $[V]$: actually, if $C\cdot\ha{H}=0$, the previous 
arguments show that either $C \subset B$ or $C$ is contained in a fiber of
$q$, both are $V$-rationally connected, hence $C$ is numerically 
proportional to $V$.  
\end{proof}
Unfortunately, $B$ is not $V$-rationally connected in general as shown
by the following example.
\begin{example}[see \cite{kachi} Example 11.1 and references therein]
\label{isolP2}
Fix a point $p_0$ in $\mathbb P^{3}$ and let
$$P_0:=\{\Pi\in(\mathbb P^{3})^*\,|\,p_0\in\Pi\}\simeq\mathbb P^{2}$$
be the variety of $2$-planes in $\mathbb P^3$ containing $p_0$.
Consider the variety $X\subset\mathbb P^{3}\times P_0$ defined as
$$ X:=\{(p,\Pi)\in\mathbb P^{3}\times P_0\,|\,p\in\Pi\}.$$
Then $X$ is a smooth Fano $4$-fold, with Picard number $2$ and
pseudo-index $2$. The two elementary extremal contractions are given by
the projections on the two factors.

The morphism $X\to P_0$ is a fibration in $\mathbb P^{2}$: the fiber over a
point is the plane corresponding to that point.

Consider the morphism $X\to \mathbb P^{3}$. If $p\neq p_0$, the fiber over
$p$ is the $\mathbb P^{1}$ of planes containing $p$ and $p_0$. But the fiber
$F_0$ over $p_0$ is naturally identified with $P_0$, hence it is
isomorphic to $\mathbb P^{2}$.
We have $\ma{N}_{F_0/X}=\Omega^1_{\mathbb P^{2}}(1)$ and
$(-K_X)_{|F}=\ol_F(2)$.
$$
\xymatrix{ {\ma{C}}\ar[d]_{\pi = p}\ar[r]^{F=e} &
  X\ar[d]^{q'}\ar@{-->}[dl]_{q}  
\\ V \ar[r]_{\psi}& 
{\mathbb P^{3}}}$$ 
Here $V\to\mathbb P^{3}$ is the blow-up of $p_0$ and $\ma{C}\to X$ is the
blow-up of $F_0$. Observe that $V$ is a family of extremal
irreducible 
rational curves of anticanonical degree $2$.

If we consider $X\times\mathbb P^{1}$ with the same family of curves, 
we have $\dim Y=4$, $f_V=1$ and
$B=F_0\times\mathbb P^{1}$ which is  not $V$-rationally connected.
\end{example}
We finally get the following result: if $B$ has the smallest
possible dimension, then it is $V$-rationally connected. 
\begin{lemma}
\label{dimT}
Let $X$ be a normal and $\Q$-factorial projective variety, and $V$ 
be a
covering and quasi-unsplit family of rational 1-cycles on $X$.

If $\dim B= f_V+1$, then
every connected component of $B$ is a $V$-equivalence class.
\end{lemma}
\begin{proof}
By Proposition \ref{baselocus},
 we know that $B$ is the union of all $V$-equivalence classes whose
 dimension is $f_V+1$. Since each of these equivalence classes must
 contain an irreducible component of $B$, they
 are in a finite number, and each is contained in a connected
 component of $B$.

So if $B_0$ is a connected component of $B$, we have
$B_0=F_1\cup\cdots\cup F_r$, where each $F_i$ is a $V$-equivalence
class. We want to show that $r=1$. 

Assume by contradiction that $r>1$. Observe that the $F_i$'s are
disjoint and $B_0$ is connected, hence at least one $F_i$ is not a
closed subset of $X$, assume it is $F_1$.

Then $F_1$ is a countable union of closed subsets. Considering the
decomposition of $B_0$ as a union of irreducible components, we find
an irreducible component $T$ of $B_0$ such that
$$T=\bigcup_{m\in\mathbb{N}} K_m$$
where each $K_m$ is a non empty proper closed subset of $T$. 
Since $T$ is an irreducible complex projective variety,
this is impossible.
\end{proof}
We then reformulate in a single result what we proved so far, and show that it
implies Theorem \ref{result}.  
\begin{prop}
\label{finally}
Let $X$ be a normal and $\Q$-factorial projective variety, and $V$ a
covering and quasi-unsplit family of rational 1-cycles on $X$.
Then:
\begin{enumerate}[$(i)$]
\item either $B = \emptyset$ or
$\dim(B) \geq f_V +1$, 
\item if $B = \emptyset$ or 
if $\dim(B) = f_V +1$  
then $\R_{\geq 0} [V]$ is a geometric extremal ray
of the Mori cone $\overline{\rm NE}(X)$.
\end{enumerate}
\end{prop}
\begin{proof}[Proof of Theorem \ref{result}]
Just notice that if $f_V\geq n-3$ and
$B$ is not empty, Proposition \ref{finally} $(i)$ gives
 $\dim B\geq f_V+1\geq n-2$, so
 $\dim B=n-2=f_V+1$. Then Proposition \ref{finally} $(ii)$ gives that 
$\R_{\geq 0} [V]$ is a geometric extremal ray
of the Mori cone $\overline{\rm NE}(X)$.
\end{proof}
%%%%%%%%%%%%%%%%%%%%%%%%%%%%%%%%%%%%%%%%%%%%%%%%%%%%%%%%%%%%%%%%%%%%%%%%%%%%%%%
\section{Existence of a geometric quotient}
\noindent Let $V$ be a covering and quasi-unsplit
family of rational 1-cycles on $X$, and assume that there exists a
geometric quotient $q'\colon X\to Y'$ for $V$. 

Observe that $q'$ has the
following property:
\emph{for any irreducible curve $C$ in $X$, $q'(C)$ is a point if and
only if $[C]$ is proportional to $[V]$}.

Conversely, we show that a morphism with the property above
is quite
close to be a geometric quotient.
\begin{prop}
\label{birational}
Let $X$ be a normal and $\Q$-factorial projective variety, and $V$ a
covering and quasi-unsplit family of rational 1-cycles on $X$.

Assume that there exists a  morphism with connected fibers
$q'\colon X\to Y'$, onto a complete and
normal algebraic variety $Y'$, such that
for any irreducible curve $C$ in $X$, $q'(C)$ is a point if and
only if $[C]$ is proportional to $[V]$.

Then there exists a birational morphism $\psi\colon Y\to Y'$ that fits
into the commutative diagram:
\begin{equation}
\label{comm}
 \xymatrix{ {{Z}}\ar[d]_{p}\ar[r]^{e} & 
X \ar@{-->}[ld]_{q}\ar[d]^{q'}
\\ {Y} \ar[r]_{\psi}& {Y'} }
\end{equation}
Moreover, if $B':=q'(B)$, we have $(q')^{-1}(B')=B$, and
$$B'\,=\,\{y\in Y'\,|\,\dim(q')^{-1}(y)>f_V\}\,=\,\{y\in
Y'\,|\,\dim\psi^{-1}(y)>0\}. $$
In particular, every fiber of $q'$ over $Y'\setminus B'$ is a
$V$-equivalence class.
\end{prop}
Observe that in example \ref{isolP2}, $\psi$ is not an isomorphism.
\begin{proof}
Let's show first of all that $(q')^{-1}(B')=B$. 

If $C\subset X$ 
is an irreducible curve contained in a fiber of $q'$, then either
$C\cap B=\emptyset$, or $C\subseteq B$.
In fact, assume that $C\cap B\neq \emptyset$.
Let
$\widehat{H}_0,\dotsc,\widehat{H}_N$ be as in the proof of
Proposition \ref{baselocus}. 
Then for any $i=0,\dotsc,N$, 
we have $C\cdot \widehat{H}_i=0$ and $C\cap
\widehat{H}_i\neq\emptyset$, 
hence $C\subseteq\widehat{H}_i$ so $C\subseteq B$.

Since $q'$ has connected fibers, we see that for every fiber $F$ of
$q'$, either $F\cap B=\emptyset$, or $F\subseteq B$. This means that 
$(q')^{-1}(q'(B))=B$.

The existence of $\psi$ as in
\eqref{comm} follows easily from the normality of $Y$ and 
the fact that $q'$ contracts all
curves in $V$, hence all $V$-equivalence classes. 
Observe that $\psi$ is surjective with connected fibers.

Let's show that $p$ contracts to a point any fiber of $q'\circ e$ over
$Y'\setminus B'$.

Let $F$ be a fiber of $q'$ over $Y'\setminus B'$, 
then we have $F\subset X\setminus B$.
Let $C\subset F$ be an irreducible curve, and
choose an irreducible curve $C'\subset X\setminus B$ 
which is a component of a cycle of
 the family $V$. Since $q'(C)$ is a point, there exists
$\lambda\in\Q_{>0}$ such that
$C\equiv\lambda C'$ in $X$.

Set $X^0:=X\setminus B$. Notice that $e$ is an isomorphism over
$X^0$, so $X^0$ can be viewed also as an open subset of $Z$; in the same
way the curves $C$ and $C'$ 
can be viewed also as a curves in $Z$.
Let's
show that $C\equiv\lambda C'$ still holds in $Z$.

Let $L\in\Pic Z$, and write $L_{|X^0}=\ol_{X^0}(D)$, $D$ a Cartier
divisor in $X^0$. Let $\ov{D}$ be the closure of $D$ in $X$ (meaning, if
$D=\sum_i a_i V_i$, that $\ov{D}=\sum_i a_i \ov{V}_i$) and let
$m\in\Z_{>0}$ be such that $m\ov{D}$ is Cartier in $X$. Then set
$M:=e^*(\ol_X(m\ov{D}))\in\Pic Z$. 

By construction, $M\otimes L^{\otimes (-m)}$ is trivial on $X^0$, so
we can write $L^{\otimes m}=M\otimes\ol_Z(G)$, where $G$ is a Cartier
divisor in $Z$ with $\Supp G\subseteq E$.

Now observe that $C\cdot G=C'\cdot G=0$, because both curves are
disjoint from $E$, and that $C\cdot M=\lambda C'\cdot M$ by the
projection formula. 
Then we have $C\cdot L=\lambda C'\cdot L$, so $C\equiv\lambda C'$
in $Z$. 

Then
$p$ must contract $C$ to a point, because $Y$ is projective. 
Since $e^{-1}(F)$ is connected,
we have shown that $p$ contracts
$e^{-1}(F)$ to a point.
Since $Y$ and $Y'$ are normal, this implies that $\psi$ is an
isomorphism over $Y'\setminus B'$.

Finally, 
let $y\in B'$ and let $F'=(q')^{-1}(y)$. Then $F'\subseteq B$, so $e$
has positive dimensional fibers on $F'$, and $\dim e^{-1}(F')>\dim F'\geq
f_V$. Since $e^{-1}(F')=p^{-1}(\psi^{-1}(y))$ and $p$ has all fibers of
dimension $f_V$, we must have $\dim
\psi^{-1}(y)>0$. 
\end{proof}
We can finally prove our results.
\begin{thm}
\label{result2}
Let $X$ be a normal and $\Q$-factorial complex projective
variety of dimension $n$, having canonical singularities.
Let $V$ be  a
covering and quasi-unsplit family of rational 1-cycles on $X$.

If $\dim B\leq f_V+1$, 
then $\R_{\geq 0} [V]$ is a geometric extremal ray
of the Mori cone $\overline{\rm NE}(X)$
and the Mori contraction of $\R_{\geq 0} [V]$, 
${\rm cont}_{[V]}\colon X\to Y'$, 
is the geometric quotient for $V$. 
\end{thm}
\begin{proof}
 If $B$ is empty, then the statement is clear.
Assume that $B$ is not empty. Then Proposition \ref{finally} and Lemma
\ref{dimT}  
yield that $\dim B=f_V+1$, every connected component of $B$ 
is a $V$-equivalence class,
and $\R_{\geq 0} [V]$
is a geometric extremal ray of $\overline{\NE}(X)$.

We have to show that $-K_X\cdot[V]>0$. 
Let $V'$ be the covering family of rational 1-cycles on $X$ 
given by Lemma \ref{families}, and
consider a resolution of singularities
$f\colon X'\to X$. 
The family $V'$ determines a covering family $V''$ of rational 1-cycles in
$X'$. If
$C_0\subset X$ is a general element of the family $V'$, then
$C':=\overline{f^{-1}(C_0\setminus\Sing(X))}$ is a general element of
$V''$, and $C_0=f_*(C')$.

Since $C_0$ is reduced and irreducible, so is $C'$. Moreover $V''$ is
covering, so $C'$ is a 
free curve in $X'$, and it has
positive anticanonical degree.

Let $m\in\Z_{>0}$ be such that $mK_X$ is Cartier.
Since $X$ has canonical singularities, we have
$$mK_{X'}=f^*(mK_X)+\sum_i a_iE_i,$$
where $E_i$ are exceptional divisors of $f$ and $a_i\in\Z_{\geq 0}$. 
Then
$$-mK_X\cdot C_0=-f^*(mK_X)\cdot C'=-mK_{X'}\cdot C'+\sum_i a_iE_i\cdot
C'>0.$$
This gives $-K_X\cdot [V']>0$ 
and thus $-K_X\cdot [V] >0$.

Since $X$ has canonical singularities, the cone theorem and the
contraction theorem hold for $X$ (see \cite[Theorems 7.38 and
7.39]{debarreUT}). Moreover, 
the extremal ray $\R_{\geq
  0}[V]$ lies in the $K_X$-negative part of the Mori cone, hence it can
be contracted.

Let ${\rm cont}_{[V]}\colon X\to Y'$ be the extremal contraction; then
$Y'$ is a 
normal, projective variety, and it is $\Q$ factorial by
\cite[Proposition 7.44]{debarreUT}.

Applying 
Proposition \ref{birational}, we see that
all fibers of ${\rm cont}_{[V]}$ over $Y'\setminus {\rm cont}_{[V]}(B)$ 
are $V$-equivalence classes. 
Since connected components of $B$ are $V$-equivalence classes,
they are exactly 
 the fibers of ${\rm cont}_{[V]}$ over ${\rm cont}_{[V]}(B)$, 
and we have the statement.
\end{proof}
Observe that
Theorem \ref{toto} is a straightforward consequence of Theorem
\ref{result2}. 

%%%%%%%%%%%%%%%%%%%%%%%%%%%%%%%%%%%%%%%%%%%%%%%%%%%%%%%%%%%%%%%%%%%%%%
\section{The toric case: proof of Theorem \ref{toric}}
\noindent 
In the case $\rho_X=1$, the statement is true for 
$q'\colon X\to\{pt\}$. In fact, using Proposition \ref{birational}, we
see that the geometric quotient $Y$ must be a point. 

Assume that $\rho_X>1$.
Recall the diagram: 
$$
\xymatrix{ {\ma{C}}\ar[d]^{\pi}\ar[r]^F & X \\ V & }$$
Recall also that if $D\subset X$ is a prime invariant
Weil divisor, there is a natural inclusion
$i_D\colon\mathcal{N}_1(D)_{\R}\hookrightarrow\mathcal{N}_1(X)_{\R}$. 

\medskip

\noindent\emph{Step 1:
let $D\subset X$ be a prime invariant
Weil divisor such that $D\cdot [V]=0$. Then there exists a covering
and quasi-unsplit family 
$V_D$ of rational 1-cycles in $D$ such that 
$\R_{\geq 0}(i_D[V_D])=\R_{\geq 0}[V]$. }

Choose an irreducible component $W$ of $F^{-1}(D)$ which dominates
$D$. Set $V'_D:=\pi(W)$, and let $\mathcal{C}'_D$ be an irreducible
component of $\pi^{-1}(V'_D)$ containing $W$. Consider the
normalization 
${\ma{C}_D}$ of 
$\ma{C}'_D$, and let $\pi_D\colon
\ma{C}_D\to {V}_D$ be the Stein factorization 
of the composite map $\ma{C}_D\to \ma{C}'_D \to V'_D$.
Finally let $F_D\colon \mathcal{C}_D\to X$ be the induced map.

For $v\in V_D$, set $G_v:=F_D(\pi_D^{-1}(v))$. Then $G_v\cap
D\neq\emptyset$,
$G_v$ is connected, and 
$G_v\cdot D=0$ because $V$ is quasi-unsplit. 
This implies $G_v\subseteq D$, 
hence $F_D(\mathcal{C}_D)\subseteq D$.
Moreover, since $W$ dominates $D$, we have $F_D(\mathcal{C}_D)= D$.

Since ${V}_D$ is normal, there is a holomorphic map
${V}_D\to {\rm Chow}(D)$. Then after replacing ${V}_D$ by its
image in ${\rm Chow}(D)$ and $\ma{C}_D$ by its image in  
$\Chow(D)\times X$, we get the desired family. 

\medskip

\noindent\emph{Step 2: there
exists an invariant prime Weil divisor having intersection zero with 
$[V]$.}

In fact, let $q\colon X\dasharrow Y$ be the rational map
associated to $V$. Since $\rho_X>1$, $Y$ is not a point. Let $D$ be
a prime divisor in $Y$ intersecting $q(X^0)$ and set
$D':=\ov{q^{-1}(D)}$. Since there are curves of the family $V$ 
disjoint from $D'$, we have
$D'\cdot[V]=0$. Moreover, $D'$ is linearly equivalent
to $\sum_i a_i D_i$, where $a_i\in\Q_{>0}$
and $D_i$ are invariant
prime Weil divisors. Hence the statement.

\medskip

\noindent\emph{Step 3:
 we prove the statement.}

Let $\Sigma_X$ be the fan of $X$ in $N\cong\Z^n$,
and let $G_X$ be
the set of primitive generators of one dimensional cones in
$\Sigma_X$. It is well known that 
$G_X$ is in bijection with the set of
invariant prime divisors of $X$; for any $x\in G_X$, we denote
$D_x$ the associated divisor. Recall that for any class $\gamma\in
\mathcal{N}_1(X)_{\Q}$, we have 
$$\sum_{x\in G_X}(\gamma\cdot
D_x) x=0\ \text{ in } N\otimes_{\Z}\Q,$$ 
and that the association $\gamma\mapsto 
\sum_{x\in G_X}(\gamma\cdot
D_x) x$ gives a canonical identification
of
$\mathcal{N}_1(X)_{\Q}$ with the $\Q$-vector space
of linear relations with rational coefficients among  $G_X$.

Let $m_1x_1+\cdots+m_hx_h=0$ be the relation corresponding to
$[V]$, with $x_i\in G_X$ and $m_i$ non zero rational numbers for all $i$.
Since $V$ is covering and quasi-unsplit, 
all $m_i$'s must be positive.
 For $y\in G_X$, we have  $D_y\cdot[V]=0$ if and only if $y$ is
 different from $x_1,\dotsc,x_h$.
So by Step 2, we know that $G_X\setminus\{x_1,\dotsc,x_h\}$ is
non empty. 

The following two
statements are equivalent (see \cite[Theorem 2.4]{reid}
and \cite[Theorem 2.2]{contr}): 
\begin{enumerate}[$(a)$]
\item 
there exists a $\Q$-factorial, projective toric variety
$Y'$, and a flat, 
equivariant morphism $q'\colon X\to Y'$, such that for any curve $C$
in $X$, $q'(C)$ is a point if and only if $[C]$ is proportional to 
$[V]$;
\item
for any $\tau\in\Sigma_X$ such that  
$x_1,\dotsc,x_h\not\in\tau$, we have
\begin{equation}\label{in}
\tau+\langle
x_1,\dotsc,\check{x}_i,\dotsc,x_h\rangle\in\Sigma_X\quad\text{ for all
}i=1,\dotsc,h.\end{equation}
\end{enumerate}
Let's show $(b)$ by induction on the dimension of $X$.

Clearly, it is enough to check \eqref{in} for any maximal $\tau$ in
$\Sigma_X$ not containing any $x_i$.
Since $\{x_1,\dotsc,x_h\}\subsetneq
G_X$, such a maximal $\tau$ 
will have positive dimension.

Let $y\in G_X\cap\tau$.
We have $D_{y}\cdot[V]=0$, so by Step 1 
there exists a 
quasi-unsplit, covering family $V_{D_{y}}$ in $D_{y}$ such that
$i_{D_y}[V_{D_y}]$ is proportional to $[V]$.

 Set $\overline{N}:=N/\Z\cdot y$ and for any $z\in N$, write
$\overline{z}$ for its image in $\overline{N}$. The fan
$\Sigma_{D_{y}}$ of $D_y$ is given by the projections in
$\ov{N}\otimes_{\Z}\Q$ of all cones of $\Sigma_X$ containing $y$. 
The relation corresponding to $[V_{D_y}]$ is
$\lambda m_1\overline{x}_1+\cdots+\lambda m_h\overline{x}_h=0$, for
some $\lambda\in\Q_{>0}$.
 By induction, we know that $(b)$ holds for $V_{D_{y}}$ in
$D_{y}$.
In particular, the projection $\ov{\tau}$ of $\tau$ 
is in $\Sigma_{D_{y}}$, so we have
$$ \ov{\tau}+\langle
\ov{x}_1,\dotsc,\check{\ov{x}}_i,\dotsc,\ov{x}_h
\rangle\in\Sigma_{D_{y}}\quad\text{ for all
}i=1,\dotsc,h.$$ 
This yields \eqref{in}.

Finally, since $q'$ is flat, 
all fibers must be $V$-equivalence classes and $B=\emptyset$.

\bigskip

\small
%\bibliographystyle{alpha}
%\bibliography{Biblio}

\bigskip

\noindent L.B. \emph{e-mail: bonavero@ujf-grenoble.fr}

\noindent S.D. \emph{e-mail: druel@ujf-grenoble.fr}

\smallskip

\noindent\emph{Institut Fourier, UFR de Math\'ematiques, Universit\'e
  de Grenoble 1, UMR 5582, BP 74, 38402 Saint Martin d'H\`eres, FRANCE}

\smallskip

\noindent C.C. \emph{e-mail: casagrande@dm.unipi.it}

\smallskip

\noindent\emph{Universit\`a di Pisa, Dipartimento di Matematica, Largo
Bruno Pontecorvo 5, 56127 Pisa, ITALY}

\end{document}